# Iteratively algebraic orders

P.H. Rodenburg
*Institute for Informatics, University of Amsterdam*

A short proof of a theorem of M.H. Albert, and its application to lattices.

## Introduction

In a preliminary version of [1], Grätzer asked if there exists a lattice that is isomorphic to its lattice of ideals and in which not every ideal is principal. This question was answered in the negative by D. Higgs [2].

A related question was posed by H.-E. Hoffmann in [3]: whether an algebraic order (poset) whose compact elements again form an algebraic order and so on, can have noncompact elements. It was answered in the negative by M.H. Albert [4].

Having forgotten a crucial element of Higgs' proof, and unable to understand the proof in [4], I invented a simpler proof, the second, and hopefully this time correct, version of which I present below, after a brief rehearsal of definitions. At the end I will indicate what baffled me.

## Substance

This note is about ordered sets, that is, sets with a (partial) ordering relation $\leq$ on them. We write $x < y$ if $x \leq y$ and $x \neq y$, and $x \prec y$ if $y$ is an *upper cover* of $x$, that is, $x < y$ and if $x \leq z \leq y$ then $z$ is either $x$ or $y$. A *chain* is a linearly ordered set.

**Definition 1**. A subset $X$ of an ordered set $L$ is *directed* if every finite subset of $X$ has an upper bound in $X$.

In particular, the void subset of a directed set has an upper bound, so directed sets are nonvoid.

**Definition 2**. An ordered set $L$ is *upwards complete* if every directed $X \subseteq L$ has a supremum in $L$.

The supremum of $X$ is denoted by $\bigvee X$. We write $(X]_L$, omitting the subscript if it can be derived from the context, for

$$\{y \in L \mid \text{for some } x \in X, y \leq x\}.$$

Instead of $(\{x\}]$, we write $(x]$. Dually we have $[X)$ and $[x)$.

**Definition 3**. An element $k$ of an ordered set $L$ is *compact* if for every directed $X \subseteq L$, $k \leq \bigvee X$ implies $k \in (X]$.

We denote the set of compact elements of an ordered set $L$ by $K(L)$. We put $K^0(L) = L$, $K^{n+1}(L) = K(K^n(L))$.

**Definition 4**. An ordered set $L$ is *algebraic* if it is upwards complete and for every $x \in L$, $(x] \cap K(L)$ is directed and $x$ is its supremum. It is *iteratively algebraic* if for all $n$, $K^n(L)$ is algebraic.

**Theorem**. If an ordered set $A$ is iteratively algebraic, $A = K(A)$.

**Proof**. Assume $A$ is iteratively algebraic, and $a \in A - K(A)$. Since the supremum of a chain of noncompact elements is noncompact, by Zorn's Lemma, $A$ contains a *maximal* noncompact element $m_0 \geq a$. Clearly, $[m_0)$ satisfies the ACC — the supremum of an infinitely ascending chain cannot be compact.

The element $m_0$ is the supremum of a directed set $C$ of compact elements. We know that $m_0$ is not in $K(A)$; but $C \subseteq K(A)$, and since $K(A)$ is algebraic, $C$ has a supremum $a'$ in $K(A)$. Now $a'$ is noncompact in $K(A)$. Again using Zorn's Lemma, we find a maximal noncompact $m_1 \geq a'$ in $K(A)$. Then $m_1 \geq a' > m_0$; and repeating the argument we get an infinitely ascending chain

$$m_0 < m_1 < m_2 < \ldots$$

in $[m_0)$, contradicting the ACC. So $A = K(A)$. ☒

**Corollary**. If a lattice $L$ is isomorphic to its ideal lattice, all its ideals are principal.

**Proof**. As an ordered set, the ideal lattice Idl($L$) is algebraic; so likewise $L$ is algebraic. ($L$ will even be an algebraic *lattice* if it has a 0.) The compact elements of Idl($L$) are the principal ideals. The sublattice of principal ideals is obviously isomorphic to $L$, so $L \cong K(L)$, which implies that $L$ is iteratively algebraic. Then by the Theorem, all the elements of $L$, and hence all the elements of Idl($L$), are compact. ☒

**Discussion**

The proof of the theorem certainly owes to Higgs, but omits his main idea: a construction of double sequences of compact elements, based on the observation that a lower cover of an ideal generated by a compact element must be principal. Albert [4] claims to prove the theorem, but his conclusion that $A = K(A)$, after a transfinite lopping off of maximal elements, appears right out of the blue. Hansoul [5] suggests a proof of the dual of Albert's theorem along the lines of [2].